\def\e{{\epsilon}}
\def\t{{\tau}}
\def\Th{{\Theta}}
\def\w{{\omega}}
\def\W{{\Omega}}

\def\E{{\rm\bf E}}
\def\F{{\cal F}}
\def\G{{\cal G}}

\def\P{{\cal P}}

\def\l({{\left(}}
\def\r){{\right)}}
\def\cpf{{$\diamond$}}
\def\pf{{$\Box$}}
\def\'{{^\prime}}

\def\({{\Biggl(}}
\def\){{\Biggr)}}
\def\[{{\Biggl[}}
\def\]{{\Biggr]}}

\documentclass[12pt]{article}
\newtheorem{theorem}{Theorem}

\newtheorem{conjecture}[theorem]{Conjecture}

\newtheorem{problem}[theorem]{Problem}
\newtheorem{proposition}[theorem]{Proposition} 
\usepackage{latexsym,psfig,multicol,amssymb}

\begin{document}

%
%
\title{On the Pebbling Threshold Spectrum}
\author{
Andrzej Czygrinow\\
Department of Mathematics and Statistics\\
Arizona State University\\
Tempe, Arizona 85287-1804\\
email: andrzej@math.la.asu.edu\\
 \\
and \\
 \mbox{}\\
Glenn Hurlbert\thanks{Partially supported by National Security Agency
grant \#MDA9040210095.}\\
Department of Mathematics and Statistics\\
Arizona State University\\
Tempe, Arizona 85287-1804\\
email: hurlbert@asu.edu\\
}
\maketitle
\newpage

%
%
\begin{abstract}
A configuration of pebbles on the vertices of a graph is solvable if one
can place a pebble on any given root vertex via a sequence of pebbling steps.
A function is a pebbling threshold for a sequence of graphs if a randomly 
chosen configuration of asymptotically more pebbles is almost surely 
solvable, while one of asymptotically fewer pebbles is almost surely not.
In this note we show that the spectrum of pebbling thresholds for graph 
sequences spans the entire range from $n^{1/2}$ to $n$. 
This answers a question of Czygrinow, Eaton, Hurlbert and Kayll.
What the spectrum looks like above $n$ remains unknown.
\vspace{0.2 in}

\vspace{0.2 in}

\end{abstract}
\newpage

%
%
\section{Introduction}\label{Intro}
Let $G= (V,E)$ be a connected graph on $n$ vertices and let $C$ be a 
{\it configuration} of $t$ unlabeled pebbles on $V$ 
(formally $C$ is multiset of $t$ elements from $V$). 
A {\it pebbling step} consists of removing two pebbles from a vertex $v$ and 
placing one pebble on a neighbor of $v$.
A configuration is called $r$-{\it solvable} if it is possible to move 
at least one pebble to vertex $r$ by a sequence of pebbling steps. 
A configuration is called {\it solvable} if it is 
$r$-solvable for every vertex $r \in V$. 
The {\it pebbling number} of $G$ is the smallest integer $\pi(G)$ such that 
every configuration of $t=\pi(G)$ pebbles on $G$ is solvable. 
Pebbling problems have a rich history and we refer to \cite{Hur} for a 
thorough discussion. 

Let us now recall some asymptotic notation that will be used in the paper.
For two functions $f=f(n)$ and $g=g(n)$, we write $f \ll g$ (or 
$f \in o(g)$) if $f/g$ approaches zero as $n$ approaches infinity, 
$f \in O(g)$ ($f \in \W(g)$) if there exist positive constants 
$c,k$ such that $f < cg$ ($f > cg$) whenever $n>k$. 
We will also use $f \sim g$ if $f/g$ approaches 1 as $n$ approaches infinity.
Finally to simplify the exposition we shall always assume, whenever needed, 
that our functions take integer values.

In this note, we will be interested in the following random model 
introduced in \cite{CEHK}. 
A configuration $C$ of $t$ pebbles assigned to  $G$ is selected 
randomly and uniformly from all ${ n+t -1 \choose t}$ configurations. 
The problem to investigate, then, is to find what values of $t$, as functions
of the number of vertices $n=n(G)$, make $C$ almost surely solvable.
More precisely, a function $t=t(n)$ is called a {\it threshold} of a 
graph sequence $\G=(G_1,\ldots,G_n,\ldots)$, where $G_n$ has $n$ vertices,
if the following conditions hold as $n$ tends to infinity:
\begin{enumerate}
\item for $t_1 \ll t$ the probability that a configuration of 
$t_1$ pebbles is solvable tends to zero, and
\item for $t_2 \gg t$ the probability that a configuration of 
$t_2$ pebbles is solvable tends to one. 
\end{enumerate}
We denote by $\t(\G)$ the set of all threshold functions of $\G$.
It is not immediately clear, however, that $\t(\G)$ is nonempty for all $\G$.
Nonetheless it is proven to be the case in \cite{BBCH}.

Note that the model defined above is different than an ``independent'' model,
in which each pebble independently selects a vertex on which to be placed. 
Indeed the difference is not merely a technical issue. 
For example, consider $P_n$, the path on $n$ vertices. 
In the independent model it is trivial to show that almost surely a 
configuration with $t \gg n\lg{n}$ pebbles is solvable. 
In this ``dependent'' model the situation is completely different. 
It is proven in \cite{BBCH} that if $t \ll n 2^{\sqrt{(\lg n)/2}}$ 
then a configuration is almost surely not solvable. 
Establishing an exact threshold for a graph in the dependent model 
is usually not a trivial task. 
In fact, even for the sequence of paths $\P=(P_1,\ldots,P_n,\ldots)$,
no exact threshold is known. 
The best results to date are 
$\t(\P)\subset\W(n2^{\sqrt{c\lg n}})\cap O(n2^{\sqrt{d\lg n}})$ for any
$c<1/2$ and any $d>1$.
The lower bound is found in \cite{BBCH} and the upper bound is found in
\cite{GJSW}.
The main purpose of this note is to investigate what functions $t=t(n)$ 
can be a pebbling threshold for some sequence of graphs.
In particular, we verify the following conjecture posed in \cite{CEHK}.  

\begin{conjecture}\label{Conj}
For every $\W(n^{1/2})\ni t_1\ll t_2\in O(n)$ there exists a graph sequence
$\G=(G_1,\ldots,G_n,\ldots)$ such that $\t(\G)\subset\W(t_1)\cap O(t_2)$.
\end{conjecture}

Let $m$ be an positive integer, $W= \{v_1,v_2, \dots v_m\}$ and, 
$S=\{v_{m+1}, \dots, v_n\}$. 
Consider the graph $F_{m,n}=(V, E)$, where the set of vertices $V=W\cup S$
and the set of edges $E$ is defined as follows:
for every $i=1, \dots, m-1$, $\{v_i,v_{i+1}\}\in E$ and 
for $i=m+1, \dots, n$, $\{v_i, v_m\} \in E$. 
In other words $F_{m,n}$ is a path on $m$ vertices with a 
star on $n-m+1$ vertices attached to one of the endpoints of the path.
(We like to think of $F$ as a {\it fuse} with {\it wick} $W$ and sparks $S$.)
Finally, for $m$ a function of $n$,
define the graph sequence $\F_m=(F_{m,1},\ldots,F_{m,n},\ldots)$.

\begin{theorem}\label{main}
Let $\e=\e(n)<1/2$ be such that $n^{1-\e} \ll n$. 
Then for $m=(1-2\e)\lg n$ we have $\t(\F_m)=\Th(n^{1-\e})$.
\end{theorem}

Note that Theorem \ref{main} implies Conjecture \ref{Conj}. 
Indeed, for given $t\in\W(t_1)\cap O(t_2)$ it is enough to consider 
$\F_m$ with $m = \lg{\frac{t^2}{n}}$.
The rest of the note is devoted to a proof of Theorem \ref{main}.

%
%
\section{Proof of Theorem \ref{main}}\label{Proof}
In this section, we prove Theorem \ref{main}. 
We divide the argument into two propositions. 
In the first one we show the upper bound, while the second 
contains a proof of the lower bound.
Let $F_{m,n} = (W\cup S,E)$ be a graph defined above. 
Assume that $ t \ll n$ and let $C$ be a configuration of $t$ pebbles on
$F$. Since $C(v_1) + C(v_2) + \dots C(v_n) = t$, we have the expectation
\begin{equation}\label{e0}
\E[C(v_i)] = \frac{t}{n}\ .
\end{equation}
First, for a fixed vertex $v$ and $i \geq 1$, we compute the
probability 
$$\Pr[C(v)=i] = \frac{{n+t-i-2 \choose t-i}}{{ n+t-1 \choose t}}\ .$$
We next compute
$${n+t-i-2 \choose t-i} = 
\[\(\frac{t-i+1}{n+t-i-1}\)\cdots\(\frac{t}{n+t-2}\)\] { n+t -2 \choose t}$$
$$= \[\(\frac{t-i+1}{n+t-i-1}\) \cdots \(\frac{t}{n+t-2}\)\] 
	\(\frac{n-1}{n+t-1}\) { n+ t-1 \choose t}\ .$$
This yields
\begin{eqnarray*}
\(\frac{n-1}{n+t-1}\) \(\frac{t-i}{n+t-i}\)^i {n+t-1 \choose t} 
	&\leq&{n+t-i-2 \choose t-i} \\
&\leq&\(\frac{t}{n}\)^i {n+t-1 \choose t}\ .
\end{eqnarray*}
Therefore,
\begin{equation}\label{e1}
\(\frac{n-1}{n+t-1}\)\(\frac{t-i}{n+t-i}\)^i 
	\leq \Pr[C(v)=i] 
	\leq \(\frac{t}{n}\)^i.
\end{equation}

\begin{proposition}\label{upper}
Let $\e= \e(n) \geq 0$ and let $\w = \w(n) \rightarrow \infty$ be such that 
$t = \w(n) n^{1-\e} \ll n$. 
Then for $m=(1-2 \e)\lg{n}$ the probability that a random configuration 
of $t$ pebbles on $F_{m,n}$ is solvable approaches one as $n$ approaches 
infinity.
\end{proposition}

\noindent
{\it Proof}.  
Let $F_{m,n} =(W\cup S,E)$, where $m=(1-2 \e)\lg{n}$, $W=\{v_1, \dots, v_m\}$ 
and $S=\{v_{m+1}, \dots , n\}$.
Let $L_2 = \{v| C(v) = 2\}$ and consider $X = |S \cap L_2|$. Then $ X =
\sum_{i=m+1}^{n} X_i$, where $X_i = 1$ if and only if $C(v_i)=2$. By (\ref{e1}),
$$ \E[X] \leq |S| \(\frac{t}{n}\)^2$$
and
$$\E[X] \geq |S|\(\frac{n-1}{n+t-1}\)\(\frac{t-2}{n+t-2}\)^2\ .$$
Since $t \ll n$, we have
\begin{equation}\label{e2}
\E[X] \sim (n - (1-2\e)\lg{n})(\w(n)n^{-\e})^2 \sim \w(n)^2 n^{1 - 2\e}\ . 
\end{equation}  
Recall that $v_m$ denotes the center of the set $S$. 
We shall show that $\Pr[X \geq n^{1-2\e}] \rightarrow 1$. 
Then we can accumulate $n^{1 - 2 \e}$ pebbles on $v_m$, 
and since $m= (1-2 \e)\lg{n}$
we can pebble from $v_m$ to any other vertex of $F_{m,n}$. 
Indeed, 
$$\sigma_X^2 = \E[X^2] - \E[X]^2 
= \sum_{i=m+1}^{n} \E[X_i^2] + \sum_{i \neq j} \E[X_iX_j] - \E[X]^2\ ,$$
and since $\E[X_iX_j] \leq \E[X_i]\E[X_j]$, 
$$\sigma_X^2 \leq \sum_{i=m+1}^{n} \E[X_i] = \E[X]\ .$$
Using (\ref{e2}), we have 
$\Pr[X < n^{1-2 \e}] \leq \Pr[\ |X-\E[X]|> \E[X]/2\ ]$,
which by Chebyshev's inequality is at most
$$\frac{4}{\E[X]} \rightarrow 0\ .$$    

\hfill \cpf

\begin{proposition}\label{lower}
Let $\e= \e(n) \geq 0$ and let $\w = \w(n) \rightarrow \infty$.  
If $t = \frac{n^{1-\e}}{\w}$ then for $m=(1-2 \e)\lg{n}$ the probability
that a random configuration of $t$ pebbles on $F_{m,n}$ 
is solvable approaches zero as $n$ approaches infinity.
\end{proposition}

\noindent
{\it Proof}.
Let $F_{m,n} =(W\cup S,E)$, where $m=(1-2 \e)\lg{n}$, $W=\{v_1, \dots, v_m\}$ 
and $S=\{v_{m+1}, \dots , v_n\}$.
Set $L_i = \{v| C(v) = i\}$ .
Then $\E[|L_i \cap S|]\leq |S| (\frac{t}{n})^i$ and so
\begin{equation}\label{e3}
\E[|L_i \cap S|] \leq \frac{n - (1 -2\e)\lg{n}}{[\w(n)n^{\e}]^i}\ . 
\end{equation}
Let $A$ be the number of pebbles that can be accumulated on $v_m$ using the
pebbles assigned to vertices from $S$. 
Then
\begin{equation}\label{e4}
\E[A] = \E[|S\cap L_2|] + \E[|S \cap L_3|] + 2\E[|S \cap L_4|] + \dots
	\lfloor \frac{t}{2} \rfloor \E[|S \cap L_t|]\ . 
\end{equation}
Using (\ref{e3}) we can bound $\E[A]$ from above by
\begin{eqnarray}
\E[A]&<&\frac{n -(1 -2\e)\lg{n}}{[\w n^{\e}]^2} 
	\sum \frac{(k+1)}{[\w(n)n^{\e}]^{k}}\label{en1a}\\
&<&\frac{2(n - (1 -2\e)\lg{n})}{[\w n^{\e}]^2}\label{en1b}\\
&<&\frac{2n^{1-2\e}}{\w^2}\ .\label{en1c}
\end{eqnarray}

Define the following random variable
$$Y = \sum_{k=0}^{m-1}\frac{C(v_{k+1})}{2^k} + \frac{A}{2^{m-1}}$$
and note that $Y \geq 1$ if and only if $C$ is $v_1$-solvable. 
Then by (\ref{e0})
$$\E[Y]  \leq \frac{2}{n^{\e}\w} + \frac{\E[A]}{2^{m-1}}\ ,$$
and  by (\ref{en1a}-\ref{en1c})
$$\E[Y] < \frac{2}{n^{\e}\w} +\frac{n^{1 -2\e}}{\w^2 2^{m-2}} =
\frac{2}{n^{\e}\w} +\frac{4}{\w^2} \rightarrow 0\ .$$
Therefore, by Markov's inequality,
$$\Pr[Y \geq 1] \leq \E[Y] \rightarrow 0\ .$$

\hfill \cpf

\noindent
{\it Proof of Theorem \ref{main}}.
By Proposition 1 and Proposition 2, for $m=(1 -2\e)\lg n$,
$$\t(\F_{m,n})=\Th(n^{1-\e}).$$
\hfill \pf

%
%
\section{Remarks}
We finish with a few open problems and conjectures.
The most obvious question remaining in this work is what functions from
$\W(n)$ can be a threshold for some graph sequence.

\begin{conjecture}
Let $t\in\t(\P)$ be a threshold for the sequence of paths.
For every $\W(n)\ni t_1\ll t_2\in O(t)$ there exists a graph sequence 
$\G=(G_1,\ldots,G_n,\ldots)$ such that $\t(\G)\subset\W(t_1)\cap O(t_2)$.
\end{conjecture}

It is shown in \cite{CEHK} that there are no thresholds $\t\gg t$
for $t\in\t(\P)$.
That is, the sequence of paths has the highest of all thresholds.
As mentioned previously, however, this threshold is unknown, lying in
the range $\W(n2^{\sqrt{c\lg n}})\cap O(n2^{\sqrt{d\lg n}})$ for any
$c<1/2$ and any $d>1$.

\begin{problem}
Find the threshold $\t(\P)$ for the sequence of paths.
\end{problem}

Finally, we mention a problem that is somewhat tangential.
In this note we investigated the almost sure solvability of random pebbling
configurations.
That is, we considered random configurations on a given graph.
One can also consider given configurations on a random graph.
In other words, what almost surely is the pebbling number of a random graph?
Along these lines, all graphs satisfy $\pi(G)\ge n(G)$, and those for which
$\pi(G)=n(G)$ are called {\it Class 0}.
It is proven in \cite{CHKT} that the random graph threshold (for the
uniform and independent probability that a given edge appears) of the
Class 0 property lies in the range $\W(\lg n/n)\cap o((n\lg n)^{1/d}/n)$
for all $d>0$.

\begin{problem}
Find the random graph threshold for the Class 0 property.
\end{problem}

%
%

%
%
\bibliographystyle{plain}

\begin{thebibliography}{99}
\bibitem{BBCH} 
A. Bekmetjev, G. Brightwell, A. Czygrinow and G.H. Hurlbert,
Thresholds for families of multisets, with an application to graph pebbling, 
Discrete Math., to appear.
\bibitem{CEHK}
A.~Czygrinow, N.~Eaton, G.~Hurlbert and P.M.~Kayll,
On pebbling threshold functions for graph sequences,
Discrete Math. {\bf 247/1-3} (2002), 93--105.
\bibitem{CHKT} 
A.~Czygrinow, G.~Hurlbert, H.~Kierstead and W.T.~Trotter,
A note on graph pebbling,
Graphs and Combin., to appear.
\bibitem{GJSW} 
A.~Godbole, M.~Jablonski, J.~Salzman and A.~Wierman,
An improved upper bound for the pebbling threshold of the $n$-path,
preprint (2002).
\bibitem{Hur} 
G.H.~Hurlbert,
A survey of graph pebbling,
Congress. Numer. {\bf 139} (1999), 41--64.

\end{thebibliography}
%

%
%
\end{document}